\documentclass[12pt]{amsart}
\usepackage{graphicx}
\usepackage{multicol}
\usepackage[normalem]{ulem}

\usepackage{amsmath,amsfonts,amsthm,amssymb,amscd}
\usepackage[mathscr]{eucal}
\allowdisplaybreaks
\usepackage{hyperref}
\usepackage{times}
\usepackage[normalem]{ulem}
\usepackage[german,scrtime]{prelim2e}

\usepackage[curve,matrix,arrow]{xy}

\newcommand{\amsv}[1]{\noindent
  \mbox{}\kern-.05\textwidth
  \framebox{\sl{\footnotesize{\begin{minipage}[t]{1.1\textwidth} #1
        \end{minipage}}}}}

\newcommand{\amsint}[1]{\noindent
  \mbox{}\kern-.08\textwidth
  {\sl{\footnotesize
      internal:}}\hbox to .3\textwidth{\hrulefill}\\*[-.4\baselineskip]
  {\begin{minipage}[t]{1.1\textwidth}#1\\*[-.8\baselineskip]
      \mbox{}\kern-.08\textwidth\hbox to .35\textwidth{\hrulefill}
      \medskip
    \end{minipage}}}

\renewcommand{\hat}{\widehat}

\newcommand{\daha}{{\mathcal H}\kern-6pt{\mathcal H}}

\newcommand{\ii}{\mathbf i}

\newcommand{\idem}{\boldsymbol{e}}
\newcommand{\nilp}{\boldsymbol{w}}

\newcommand{\fidem}{\boldsymbol{f}}
\newcommand{\unilp}{\boldsymbol{u}}

\renewcommand{\leq}{\,{\leqslant}\,}

\renewcommand{\le}{\,{\leqslant}\,}

\newcommand{\q}{\mathfrak{q}}
\newcommand{\tc}{\mathfrak{t}}

\newcommand{\fusion}{%
  \mathop{{\otimes}\kern-7pt\raisebox{.6pt}{%
      \mbox{\footnotesize${\bullet}$}}}}

\newcommand{\UresSL}[1]{\overline{\mathscr{U}}_{\q} s\ell(#1)}

\newcommand{\ffrac}[2]{\mbox{\footnotesize$\displaystyle\frac{#1}{#2}$}}
\newcommand{\half}{%
  \mathchoice{\ffrac{1}{2}}{\frac{1}{2}}{\frac{1}{2}}{\frac{1}{2}}}

\newcommand{\cchi}{\boldsymbol{\chi}}

\newcommand{\vvarphi}{\hat{\pmb{\varphi}}}

\newcommand{\radmap}{\pmb{\boldsymbol{\phi}}}

\newcommand{\fdec}{{\rm f}}
\newcommand{\udec}{{\rm u}}
\newcommand{\kdec}{{\rm k}}

\newcommand{\modS}{\mathscr{S}}

\newcommand{\rep}{\mathscr}  

\newcommand{\oC}{\mathbb{C}}

\newcommand{\oZ}{\mathbb{Z}}

\newcommand{\ribbon}{{\boldsymbol{v}}}

\newcommand{\Ss}{\mathop{\mathsf{s}}}
\newcommand{\Pp}{\mathop{\mathsf{p}}}

\newcommand{\cheb}{{U}}

\numberwithin{equation}{section} \makeatletter
\@addtoreset{equation}{section}
\@addtoreset{subsubsection}{section}

\def\@secnumfont{\bfseries}
\def\subsubsection{\@startsection{subsubsection}{3}%
  \z@{.5\linespacing\@plus.7\linespacing}{-.5em}%
  {\normalfont\bfseries}}
\def\paragraph{\@startsection{paragraph}{4}%
  \z@\z@{-\fontdimen2\font}%
  \normalfont\bfseries}
\def\subparagraph{\@startsection{subparagraph}{5}%
  \z@\z@{-\fontdimen2\font}%
  \normalfont\bfseries}

\makeatother

\swapnumbers
\newtheorem{Thm}[subsection]{Theorem}

\newtheorem{lemma}[subsubsection]{Lemma}

\newtheorem{prop}[subsubsection]{Proposition}

\theoremstyle{definition}

\begin{document}

\title[DAHA and LCFT]{%
  \vspace*{-4\baselineskip}
  \mbox{}\hfill\texttt{\small\lowercase{hep-th}/\lowercase{yymmxxx}}
  \\[\baselineskip]
  Double affine Hecke algebra in logarithmic conformal field theory}

\author[Mutafyan]{G.~Mutafyan}%

\address{\mbox{}\kern-\parindent gm: Moscow Institute of Physics
and Technology
 \hfill\mbox{}\linebreak \texttt{georg21a@yahoo.com}}

\author[Tipunin]{I.Yu.~Tipunin}
\address{\mbox{}\kern-\parindent iyt: Lebedev Physics
  Institute \hfill\mbox{}\linebreak \texttt{tipunin@lpi.ru}}

\begin{abstract}
  We construct the representation of Double Affine Hecke Algebra whose
  symmetrization gives the center of the quantum group $\UresSL2$ and
  by Kazhdan--Lusztig duality the Verlinde algebra of $(1,p)$ models
  of logarithmic conformal field theory.
\end{abstract}

\maketitle

\section{Introduction}\label{introd}
Recently, quantum group methods led (see the recent
review~\cite{[Semikhatov]}) to a progress in logarithmic
conformal field theory~\cite{[Gurarie]}.
For the $(1,p)$ models~\cite{[GK2]}, an equivalence between
representation categories of the chiral algebra and the quantum group
$\UresSL2$ was established~\cite{[FGST-equ]} in the general framework
of the Kazhdan--Lusztig duality~\cite{[KL]}.  Remarkably, the KL
duality extends to an isomorphism between modular group
representations on the quantum group center and on the space of
generalized characters~\cite{[FGST-mod]} of a $(1,p)$ model.
Moreover, the Verlinde algebra of $(1,p)$
models~\cite{[GaberdielKausch],[FHST]} (see also~\cite{[Flohr-Knuth]})
coincides~\cite{[FGST-mod]} with the Grothendiec ring of $\UresSL2$.

An unusual property of logarithmic conformal field theory is the
nonsemisimplicity of the Verlinde algebra.  However, this phenomenon
does not look extraordinary in the Double affine Hecke algebra
representation framework~\cite{[Cherednik-book]} of the Verlinde
algebra classification. It leads to a natural
conjecture~\cite{[Cherednik-privat]} that the $(1,p)$ model Verlinde
algebra can be realized in terms of a DAHA representation.  Indeed,
the representation of DAHA whose symmetrization gives the center of
the quantum group $\UresSL2$ and therefore the Verlinde algebra of
$(1,p)$ models is identified in the present paper.

\subsection{DAHA}
We consider the symplest DAHA~\cite{[Cherednik-book]} generated by $X$,
$Y$, and $T$ with the relations
\begin{gather}
\label{dha1}  TXT=X^{-1},\quad TY^{-1}T=Y,\quad XY=\q YXT^2,\\
\label{dha2}  (T-\tc^{\half})(T+\tc^{-\half})=0,\quad \tc=\q^{2}.
\end{gather}
In the paper we fix the deformation parameter
\begin{equation}
  \q=e^\frac{\ii \pi}{p},\qquad\q^\frac{1}{2}=e^\frac{\ii \pi}{2p},
\end{equation}
where $p=3,4,5,6,\dots$.  We let $\daha$ denote this algebra.  The
group $PSL(2,\oZ)$ acts by automorphisms on $\daha$
\begin{gather}\label{taupl}
\tau_+:\quad Y\rightarrow \q^{-1/2}XY,\quad X\rightarrow X,
\quad T\rightarrow T\\
\tau_-:\quad X \rightarrow \q^{1/2}YX,\quad Y\rightarrow
Y,\quad T\rightarrow T
\end{gather}
where
\begin{equation*}
 \left(
\begin{array}{cc}
1&1\\
0&1
\end{array} \right)\rightarrow \tau_+,\qquad
\left(
\begin{array}{cc}
1&0\\
1&1
\end{array} \right)\rightarrow \tau_-.
\end{equation*}
We note that the Fourier transform is given by
 \begin{gather}
   \label{fouriertrans} \sigma : \quad X \rightarrow Y^{-1}, \quad Y \rightarrow XT^2,\quad
   T\rightarrow T,\\
     \sigma=\tau_+\tau_-^{-1}\tau_+=\tau_-^{-1}\tau_+\tau_-^{-1}.\notag
 \end{gather}

\subsection{The representation}
We consider a $6p-4$-dimensional reducible but indecomposable
representation $\rep{Z}$ of~$\daha$. The representation $\rep{Z}$
contains the maximal subrepresentation $V^{-2}$, which in notations of
\cite{[Cherednik-book]} is defined as the quotient
$V^{-2}=\rep{P}/(X^{2p}+X^{-2p}-2)$, where $\rep{P}=\oC[X,X^{-1}]$ is
the standard representation of $\daha$ in the Lourent polynomials.
The $2p-4$-dimensional irreducible quotient $\rep{M}=\rep{Z}/V^{-2}$
is isomorphic to the representation $V_{2p-4}$ from
\cite{[Cherednik-book]} given by the quotient
$\rep{P}/\varepsilon_{-p+2}$, where
${\varepsilon_{-p+2}=\prod\limits_{j=2}^{p-1}(\q^{-j}X-\q^jX^{-1})}$.
We note also that $V^{-2}$ is also reducible and contains the maximal
$2p+4$-dimensional irreducible subrepresentation $\rep{W}$ and the
quotient $\rep{E}=V^{-2}/\rep{W}$ is isomorphic to $\rep{M}$.

In Sec.~\ref{sec:repZ}, we describe the representation $\rep{Z}$ by
the explicit action of operators $T$, $X$ and $Y$ in a basis. Then we
describe its structure and explicitly find the subrepresentation and
quotients.

\medskip

The representation $\rep{Z}$ bears a commutative associative
multiplication, which is described in Sec.~\ref{strucZ}. The
multiplication gives further the multiplication in the Verlinde
algebra.

The $PSL(2,\oZ)$ generators $\sigma$ and $\tau_+$ are realized as
a conjugation with some operators $\modS$ and $\ribbon$ respectively,
acting in $\rep{Z}$.  The operator $\ribbon$ acts by a multiplication
with an element from $\rep{Z}$, which is abusing notation denoted by
the same symbol~$\ribbon$.  By anology with~\cite{[Cherednik-book]} we
call $\ribbon$ the Gaussian element.

\subsection{Symmetrization} The operator $T$ has two
different eigenvalues $\q$ and $-\q^{-1}$ in $\rep{Z}$.  The
eigenspace of $T$ with the eigenvalue $\q$ is $3p-1$-dimensional.  We
let $\rep{T_\q}$ denote this eigenspace. In accordance with the
general theory~\cite{[Cherednik-book]}, $\rep{T_\q}$ is an associative
commutative algebra with multiplication induced by the multiplication
in~$\rep{Z}$ and at the same time a representation of~$SL(2,\oZ)$
induced by the~$PSL(2,\oZ)$-action in~$\rep{Z}$. The operators
$\modS$, $\ribbon$, $C=-(X+X^{-1})$ and $H=-(Y+Y^{-1})$ have well
defined restrictions on~$\rep{T_\q}$.  Now we are ready to formulate
the main result of the paper (the reader can find all needed quantum
group definitions in~\cite{[FGST-mod]}).
\begin{Thm}\label{Thm:main}
  \begin{itemize}
  \item $\rep{T_\q}$ is isomorphic to the center of $\UresSL2$ as
    associative commutative algebra and as $SL(2,\oZ)$ representation.
  \item Under the isomorphism the eigenvectors of $C$ correspond to
    Radford images and eigenvectors of $H$ correspond to Drinfeld
    images of the characters of $\UresSL2$ irreducible reprsentations.
  \item The Gaussian element $\ribbon$ coincides with the ribbon
    element of $\UresSL2$.
\end{itemize}
\end{Thm}
In Sec.~\ref{sec:symmetrization}, we describe the subspace of $T$ with
the eigenvalue $\q$, and in Sec. ~\ref{Proof:main} we give the proof of
theorem \ref{Thm:main}. The notations in this part correspond to the
same notations in~\cite{[FGST-mod]}.

\subsection{Structure of $\rep{Z}$\label{strucZ}}
The very important information about the representation $\rep{Z}$ is
encoded in the spectra of operators $X$ and $Y$. These operators are
not diagonalizable but both have Jordan blocks of dimension $2$.  In
order to describe their Jordan structure we introduce two basises in
which operators $X$ and $Y^{-1}$ have a Jordan form. We call the first
basis the $X$-basis and the second one the $Y$-basis. Jordan forms of $X$ and
$Y^{-1}$ coincide in $\rep{Z}$.

\subsubsection{$X$-basis} The representation $\rep{Z}$
 has the basis
 \begin{equation}\label{basis-I}
   w_1\dots
w_{2p}, e_1, e_p, e_{p+1}, e_{2p};\quad  e_2\dots e_{p-1}, e_{p+2}
\dots e_{2p-1};\quad m_2\dots m_{p-1}, m_{p+2}\dots m_{2p-1}.
 \end{equation}
 The subrepresentation $\rep{W}$ is spanned by elements $w_1 \dots
 w_{2p}$, $e_1$, $e_p$, $e_{p+1}$ and $e_{2p}$.  The elements $e_2
 \dots e_{p-1}$, $e_{p+2} \dots e_{2p-1}$ ($m_2 \dots m_{p-1}$,
 $m_{p+2} \dots m_{2p-1}$) give a basis in $\rep{E}$ (in $\rep{M}$)
 under the canonical projection.  In basis \eqref{basis-I} we have
 \begin{equation}
   Xw_s=\q^sw_s, \quad Xe_s=\q^s(e_s+w_s), \quad Xm_s=\q^sm_s.
 \end{equation}
We call \eqref{basis-I} the $X$-basis.

\subsubsection{The multiplication in $\rep{Z}$}
The representation $\rep{Z}$ is endowed with a commutative associative
multiplication, which is naturally written in basis~\eqref{basis-I} as
\begin{equation}\label{multipl-I}
e_ie_j=\delta_{i,j}e_j,\quad e_iw_j=\delta_{i,j}w_j, \quad
e_im_j=\delta_{i,j}m_j, \quad
w_iw_j=w_im_j=m_im_j=0.
\end{equation}

\subsubsection{$Y$-basis}
The representation $\rep{Z}$ contains the basis
\begin{equation}\label{basis-II}
  u_1\dots u_{2p}, f_1, f_p, f_{p+1}, f_{2p},;\quad  f_2\dots f_{p-1},
f_{p+2} \dots f_{2p-1};\quad k_2\dots k_{p-1}, k_{p+2}\dots k_{2p-1}
\end{equation}
in which $Y^{-1}$ acts as follows
\begin{equation}\label{Y-jordanbasis}
Y^{-1}u_s=\q^{s}u_s,\quad
Y^{-1}f_s=\q^{s}(f_s+u_s),\quad Y^{-1}k_s=\q^{s}k_s.
\end{equation}
We call \eqref{basis-II} the $Y$-basis.

In Subsec.~\ref{sec:ybas}, we find the $Y$-basis and give
decompositions of elements from the $Y$-basis in the $X$-basis.

\subsubsection{$PSL(2,\oZ)$-action} 
The operator $\modS$ maps the $X$-basis to the $Y$-basis
\begin{equation}\label{S-action}
  \modS w_s=u_s,\quad \modS e_s=f_s,\quad \modS m_s=k_s.
\end{equation} 
In Subsec.~\ref{sec:Strans}, we establish properties of the
$\modS$-operator. By a direct calculation, using the decompositions of
the $Y$-basis in the $X$-basis, we check that this operator satisfies
all relations~\eqref{fouriertrans} and $\modS^2=\q T^{-1}$.

In terms of the $X$-basis, the Gaussian element is
\begin{equation}\label{ribbon}
 \ribbon=\sum_{s=1}^{2p}\q^{-\frac{1}{2}(s^2-1)} e_s-w_1
  +\q^{-\frac{p^2}{2}}w_{p+1}+\left(\sum_{s=2}^{p-1}
   +\sum_{s=p+2}^{2p-1}\right)\q^{-\frac{1}{2}(s^2-1)}((p-s) w_s+p m_s).
\end{equation}
The properties of this element are described in Subsec.~\ref{sec:tau}.

In the end of Sec.~\ref{sec:psl2z}, we prove the $PSL(2,\oZ)$ relations.

%
\subsection{Notation}
We introduce Chebyshov polynomials
\begin{equation}
  \cheb_s(x)=x^{s-1}+x^{s-3}+\dots+x^{-(s-3)}+x^{-(s-1)}.
\end{equation}
In what follows we often use the numbers
\begin{gather}
\{s\}=\frac{\q^s+\q^{-s}}{\q-\q^{-1}},\qquad [s]=\frac{\q^s-\q^{-s}}{\q-\q^{-1}},\\
\omega_s=\frac{p\sqrt{2p}}{[s]^2}(-1)^{p+s+1},\qquad
  \xi_s\equiv\frac{-(-1)^{p-s}p\sqrt{2p}}{\q^s-\q^{-s}},\\
[s,j]\equiv\left\{\begin{array}{ll}s,&j=0,2p,\\
(-1)^{s-1}s,&j=p,\\
\frac{[sj]}{[j]},&j \mod p\ne 0,\end{array}\right.\qquad
\{s,j\}\equiv\left\{\begin{array}{ll}0,& j \mod p =0,\\
\frac{\{sj\}}{[j]}&otherwise. \end{array}\right.
\end{gather}

\section{Representation\label{sec:repZ}}
In this section, we recall the representation $V^{-2}$
\cite{[Cherednik-book]} and then define the representation~$\rep{Z}$,
which is an extension of $V^{-2}$. Then we find a Jordan basis for $Y$ in
which $Y^{-1}$ acts by~\eqref{Y-jordanbasis}.
\subsection{Polynomial representation $V^{-2}$\label{V^{-2}}}
The representation $\rep{Z³}$ is an extension of the representation
$V^{-2}$ from \cite{[Cherednik-book]}. To describe $V^{-2}$ we recall
the standard representation \cite{[Cherednik-book]} of $\daha$ in the
space of Laurent polynomials $\rep{P}=\oC[X,X^{-1}]$.  The $\daha$
generators act as follows
\begin{gather}
  T\rightarrow\tc^{\half} \Ss+\frac{\tc^{\half}-\tc^{-\half}}{X^2-1}(\Ss-1),\quad \tc=\q^k,\\
  Y\rightarrow -\Ss\Pp T,
\end{gather}
where
\begin{equation}
  \Ss f(X)=f(X^{-1}),\quad \Pp f(X)=f(\q X)
\end{equation}
and $X$, $X^{-1}$ act by multiplication. (We note that these
formulas differ from \cite{[Cherednik-book]} by the sign in the
definition of $Y\rightarrow \Ss\Pp T$.) The representation
$V^{-2}$ is the $4p$-dimensional representation in the quotient
space $\rep{P}/(X^{2p}+X^{-2p}-2)$.

\begin{prop}\label{prop:Y}
  \begin{itemize}
  \item The operators $X$ and $Y$ have in $V^{-2}$ eigenvalues $\q^s,
    s=1\dots 2p$, each with multiplicity~2.
  \item The Jordan basis of $X$ contains functions $e_s$ and $w_s$ for
    $s=1\dots 2p$.
  \item The Jordan basis of $Y$ contains functions $u_s$ for $s=1\dots
    2p$ and $k_s$ for $s=2\dots p-1, p+2\dots 2p-1$, and functions
    $f_1$, $f_p$, $f_{p+1}$, $f_{2p}$.
  \end{itemize}
  The action of $X$ and $Y^{-1}$ in these basises is
\begin{gather}
 \label{X-actioninV}  Xw_s=\q^sw_s,{}\quad{}Xe_s=\q^s(e_s+w_s),\\
\label{Y-actioninV}Y^{-1}u_s=\q^{s}u_s,\quad
Y^{-1}f_s=\q^{s}(f_s+u_s),\quad Y^{-1}k_s=\q^{s}k_s.
\end{gather}
\end{prop}
\begin{proof}
 To describe the spectra of these
operators, we introduce functions
\begin{equation}
\begin{split}
   w_s&=\frac{1}{4p^2}(X^{2p}-1)\sum_{j=0}^{2p-1}\q^{-sj}X^j,\\
   e_s&=\frac{1}{2p}+\frac{1}{4p^2}\sum_{j=1}^{2p-1}(2p-j)(\q^{-sj}X^j+\q^{sj}X^{-j}),
\end{split}\quad
   s=1\dots 2p,
\end{equation}
and
\begin{gather}
\label{Yeigenfunctions-u(x)}
\begin{split}
u_s&=\frac{(-1)^s}{p\sqrt{2p}}\left(\q^{s}\frac{\cheb_{p-s}(X)
   +\cheb_{p+s}(X)}{2}+\q\frac{\cheb_{p-s}(\q^{-1}X)+\cheb_{p+s}(\q^{-1}X)}{2}\right),\\
u_{p+s}&=\frac{(-1)^{p+s}}{p\sqrt{2p}}\left(\q^{p+s}\frac{\cheb_{s}(X)
     +\cheb_{2p-s}(X)}{2}+\q\frac{\cheb_{s}(\q^{-1}X)+\cheb_{2p-s}(\q^{-1}X)}{2}\right),
  \end{split}
\quad s=1\dots p,\\
\label{Yeigenfunctions-k(x)}
\begin{split}
k_s&=\frac{(-1)^{s+1}}{p\sqrt{2p}}\left(\q^{s}\cheb_{p-s}(X)+\q\cheb_{p-s}(q^{-1}X)\right),\\
k_{p+s}&=\frac{(-1)^{p+s}}{p\sqrt{2p}}\left(\q^{p+s}\cheb_{s}(X)+\q\cheb_{s}(q^{-1}X)\right),
\end{split}
\quad s=2\dots p-1\\
\label{Yeigenfunctions-f(x)} f_s=\frac{1}{p\sqrt{2p}}\cdot
\left\{\begin{aligned}
&\frac{(-1)^{p+1}\cheb_{2p}(X)}{2}, &s=p,\\
&\cheb_{p}(X),& s=2p,\\
&\frac{(-1)^p\q X\cheb_{2p}(X)}{2},&s=p+1,\\
&-q X\cheb_{p}(X),&s=1.
\end{aligned}\right.
\end{gather}
Then, \eqref{X-actioninV} is checked by a direct calculation.  It is
easy to see that $1=\sum\limits_{s=1}^{2p}e_s$. Together with
\eqref{X-actioninV} this gives
$$X^j=\sum_{s=1}^{2p}\q^{sj}(e_s+jw_s), j=0,\pm 1,\pm2\dots$$
i.e.{} functions $e_s$ and $w_s$ are linearly independent and form
a basis in $V^{-2}$.

\eqref{Y-actioninV} is also checked by a direct calculation using
the following relations
\begin{gather*}
Y\cheb_s(q^{-1}X)=(\q^s+\q^{-s})\cheb_s(q^{-1}X)-q^{-1}\cheb_s(X),\\
Y\cheb_s(X)=\q \cheb_s(q^{-1}X).
\end{gather*}
The linear independence of these vectors is proved by the standard
technic (See definition~$2.5.4$ and theorems $2.5.9$ and $2.9.3$ from
\cite{[Cherednik-book]}.).
\end{proof}

The representation $V^{-2}$ is reducible. It has a $2p+4$-dimensional
subrepresenation~$\rep{W}$ spanned by functions $w_1\dots w_{2p}$,
$e_1$, $e_p$, $e_{p+1}$ and $e_{2p}$. The quotient
$\rep{E}=V^{-2}/\rep{W}$ is isomorphic to $V_{2p-4}$ from
\cite{[Cherednik-book]}. The representation $V_{2p-4}$ is defined
in~\cite{[Cherednik-book]} as the quotient
$V_{2p-4}=\rep{P}/\varepsilon_{-p+2}$, where
${\varepsilon_{-p+2}=\prod\limits_{j=2}^{p-1}(\q^{-j}X-\q^jX^{-1})}$.

Decomposition of any polynomial $f(X)$ in the basis $e_s$, $w_s$ is
given by
\begin{gather}
\label{decomposition} f(X)=\sum_{s=1}^{2p}
\left(f(\q^s)e_s+\left(\left.X\frac{df(X)}{dX}\right|_{X=\q^s}\right)w_s\right).
\end{gather}
Using it, we check that $u_1\dots u_{2p}$, $f_1$, $f_p$, $f_{p+1}$ and
$f_{2p}$ belong to $\rep{W}$ and therefore in~$\rep{W}$ Jordan forms
of $X$ and $Y^{-1}$ coincide. But in the whole $V^{-2}$ they do not
coincide, hence automorphism~\eqref{fouriertrans} cannot be realized
as a conjugation. To recover this, we find an extension of~$V^{-2}$ to
a $6p-4$-dimensional representation $\rep{Z}$ by adding vectors
$m_2,\dots,m_{p-1}$ and $m_{p+2},\dots,m_{2p-1}$. The action of $X$ on
them is $X m_s=\q^s m_s$. The whole $\rep{Z}$ cannot be realized in a
space of polynomials in 1 varible. We describe $\rep{Z}$ in terms of
an abstract vector space.
\subsection{The representation $\rep{Z}$ in the $X$-basis}
We assume the following definition of $\rep{Z}$. The representation
$\rep{Z}$ is a $6p-4$-dimensional vector space with the basis consisting
of $4p$ vectors $e_s$ and $w_s$ with $s=1\dots 2p$, and $2p-4$ vectors
$m_s$ with $s=2\dots p-1, p+2\dots 2p-1$. The action of
$\daha$-operators in this basis is defined by the formulas:
\begin{gather}
\label{X-jordanaction}Xw_s=\q^sw_s, \quad Xe_s=\q^s(e_s+w_s), \quad Xm_s=\q^sm_s,\\
\label{Tact-w_p}Tw_p=-q^{-1}w_p-(q-q^{-1})e_p,\qquad Tw_{2p}=-q^{-1}w_{2p}-(q-q^{-1})e_{2p},\\
Tw_s=-\frac{\q^{-s}}{[s]}w_s-\frac{[s-1]}{[s]}w_{2p{-s}},\quad s\ne 0,p,\\
Te_p=\q e_p,\qquad Te_{2p}=\q e_{2p},\\
\label{Tact-e_s}Te_s=-\frac{\q^{-s}}{[s]}e_s+\frac{[s-1]}{[s]}e_{2p-s}
      +\frac{2}{(q-q^{-1})[s]^2}(w_s-w_{2p-s}), \quad s\ne 0,p,\\
Tm_{2p-1}=\q m_{2p-1}-(\q+\q^{-1})w_{1},\\
Tm_{p-1}=\q m_{p-1},\\
Tm_s=\frac{-\q^{-s}}{[s]}m_s+\frac{[s-1]}{[s]}m_{2p-s}\quad s=2\dots p-2, p+2\dots 2p-2,\\
\label{Y-w-act}
Yw_p=-q^{-1}w_{p+1}+(q-q^{-1})e_{p+1},\qquad Yw_{2p}=-q^{-1}w_1+(q-q^{-1})e_1,\\
Yw_s=-\frac{\q^{-s}}{[s]}w_{2p-s+1}-\frac{[s-1]}{[s]}w_{s+1},\quad s\ne 0,p,\\
Ye_p=-\q e_{p+1},\qquad Ye_{2p}=-\q e_1,\\
\label{Y-e-act}Ye_s=\frac{\q^{-s}}{[s]}e_{2p-s+1}-\frac{[s-1]}{[s]}e_{s+1}
 -\frac{2}{(q-q^{-1})[s]^2}(w_{s+1}-w_{2p-s+1}), \quad s\ne 0,p,\\
Ym_{2p-1}=-\q m_2-(\q+\q^{-1})w_{0},\qquad Ym_{p-1}=-\q m_{p+2},\\
\label{Y-m-act}Ym_s=-\frac{[s-1]}{[s]}m_{s+1}+\frac{\q^{-s}}{[s]}m_{2p-s+1},\quad
s=2\dots p-2, p+2\dots 2p-2.
\end{gather}
We note that the basis  $e$, $w$, $m$ by definition is the $X$-basis~\eqref{basis-I}
and~\eqref{X-jordanaction} gives the Jordan structure of $X$.

\begin{lemma}
  Operators $X$, $Y$ and $T$ defined by~\eqref{X-jordanaction}--\eqref{Y-m-act}
 satisfy the DAHA relations~\eqref{dha1} and~\eqref{dha2}.
\end{lemma}
\begin{proof}
  A direct calculation.
\end{proof}

We define a commutative associative multiplication in $\rep{Z}$ by
formulas~\eqref{multipl-I}.

\begin{prop} $\rep{Z}$ is reducible. The $2p+4$-dimensional subspace
\begin{equation}
\rep{W}\equiv\{w_1\dots w_{2p}, e_1,e_p,e_{p+1},e_{2p}\}
\end{equation}
 is invariant under the $\daha$-action and is therefore a
subrepresentation. The quotient is a direct sum:
$\rep{Z}/\rep{W}=\rep{E}\oplus\rep{M}$, where
$\rep{E}\equiv{\{e_2\dots e_{p-1},e_{p+2}\dots e_{2p-1}\}}$ and
$\rep{M}\equiv\{m_2\dots m_{p-1},$ $m_{p+2}\dots m_{2p-1}\}$.
\end{prop}
\begin{proof}
Immediately follows from \eqref{X-jordanaction}--\eqref{Y-m-act}.
\end{proof}

\subsection{$Y$-basis\label{sec:ybas}}
In this subsection we prove that the Jordan form of $Y$
is~\eqref{Y-jordanbasis}.
\begin{prop}
  A Jordan basis of $Y$ consists of $6p-4$ vectors: $4p$ vectors $f_s$,
  $u_s$ for $s=1,\dots,2p$, and $2p-4$ vectors $k_s$ for
  $s=2,\dots,p-1, p+2\dots 2p-1$. The action of $Y^{-1}$ on this
  vectors is given by~\eqref{Y-jordanbasis}.
\end{prop}
\begin{proof}
We define in the $X$-basis the vectors
\begin{equation}\label{Yeigenvectors-u}
u_s=\sum\limits_{j=1}^{2p}\udec_{j,s}^{(w)}w_j
   +\sum\limits_{j=1}^{2p}\udec_{j,s}^{(e)}e_j,\quad{}s=1\dots2p,
\end{equation}
where coefficients are
\begin{equation}
\label{u-coefficients}
\begin{split}
\udec_{j,s}^{(w)}&=\frac{(-1)^{s+j}}{\sqrt{2p}}
\Big(\q^s\{s,j\}-\q\{s,j-1\}\Big),\quad j=1\dots 2p;\\
\udec_{1,s}^{(e)}&=(-1)^s\frac{\q}{\sqrt{2p}};\quad
\udec_{2p,s}^{(e)}=(-1)^{s}\frac{\q^s}{\sqrt{2p}};\\
\udec_{p+1,s}^{(e)}&=(-1)^{p+1}\frac{\q}{\sqrt{2p}};\quad
\udec_{p,s}^{(e)}=(-1)^{p+1}\frac{\q^s}{\sqrt{2p}};\\
\udec_{j,s}^{(e)}&=0,\quad{}j\ne1,p,p+1,2p,
\end{split}
\end{equation}
the vectors
\begin{equation}
\label{Yeigenvectors-k}
k_s=\sum\limits_{j=1}^{2p}\kdec_{j,s}^{(w)}w_j+\sum\limits_{j=1}^{2p}
    \kdec_{j,s}^{(e)}e_j,\quad{}s=2\dots{}p-1,p+2\dots{}2p-1,
\end{equation}
where coefficients are
\begin{equation}
\label{k-coefficients}
\begin{split}
 \kdec_{j,s}^{(w)}&=-\frac{(p-s)}{p}\udec_{j,s}^{(w)}-\frac{(-1)^{s+j}}{p\sqrt{2p}}
        \Big(\q^s[s,j]\{1,j\}-\q[s,j-1]\{1,j-1\}\Big),\\
 \kdec_{1,s}^{(e)}&=(-1)^{s+1}{\frac{(q^s[s]+q(p-s))}{p\sqrt{2p}}},
       \quad\kdec_{p,s}^{(e)}=(-1)^{p}\frac{(q[s]+q^s(p-s))}{p\sqrt{2p}},\\
\kdec_{p+1,s}^{(e)}&=(-1)^{p}\frac{(q^s[s]+q(p-s))}{p\sqrt{2p}},
    \quad \kdec_{2p,s}^{(e)}=(-1)^{s+1}\frac{(q[s]+q^s(p-s))}{p\sqrt{2p}}, \\
\kdec_{j,s}^{(e)}&=\frac{(-1)^{s+j}}{p\sqrt{2p}}\Big(\q^s[s,j]-\q[s,j-1]\Big)\quad{}j\ne1,p,p+1,2p,
\end{split}
\end{equation}
and the vectors
\begin{equation}
\label{Yeigenvectors-f}
f_s=\sum\limits_{j=1}^{2p}\fdec_{j,s}^{(w)}w_j+\sum\limits_{j=1}^{2p}
       \fdec_{j,s}^{(e)}e_j+\left(\sum_{j=2}^{p-1}+\sum_{j=p+2}^{2p-1}\right)
                            \fdec_{j,s}^{(m)}m_j,\quad{}s=1\dots2p,
\end{equation}
where coefficients are
\begin{equation}
\label{f-coefficients}
\begin{split}
\fdec_{1,s}^{(w)}&=\frac{2(-1)^{s+1}\q^{2s}}{(\q-\q^{-1})\sqrt{2p}},\quad
\fdec_{p,s}^{(w)}=\frac{\q(-1)^{p+1}[s]}{\sqrt{2p}},\quad
\fdec_{p+1,s}^{(w)}=\frac{2(-1)^{p}\q^{2s}}{(\q-\q^{-1})\sqrt{2p}},\\
\fdec_{2p,s}^{(w)}&=\frac{\q(-1)^s[s]}{\sqrt{2p}},\\
\fdec_{j,s}^{(w)}&=-p(p-j)\kdec_{j,s}^{(e)}+\frac{(-1)^{s+j}}{\sqrt{2p}}
                        \Big(\q[s,j-1]+\q^s\{s,j\}\Big),\quad j\ne 1,p,p+1,2p,\\
 \fdec_{p,s}^{(e)}&=(-1)^{p+s+1}\fdec_{2p,s}^{(e)}=(-1)^{p+1}\frac{\q^s}{\sqrt{2p}},
             \quad \fdec_{j,s}^{(e)}=0\quad j\ne p,2p,\\
\fdec_{j,s}^{(m)}&=-p^2\kdec_{j,s}^{(e)},\quad j\ne 1,p,p+1,2p
\end{split}
\end{equation}
and the coefficient $\fdec_{j,s}^{(m)}$ in~\eqref{Yeigenvectors-f} is
$0$ for $s=1,p,p+1,2p$.

Then, \eqref{Y-jordanbasis} is checked by a simple calculation using
formulas~\eqref{Y-w-act}--\eqref{Y-m-act}.


The linear independence of these vectors is proved in the following
way. From the decompositions in the $X$-basis, we obtain that vectors
$u_s, f_1,f_p,f_{p+1}, f_{2p}$ belong to~$\rep{W}$, vectors $k_s$
belong to $\rep{W}+\rep{E}$, and vectors $f_s$ with $s=2\dots p-1$,
$p+2\dots 2p-1$ belong to $\rep{W}+\rep{M}$. We recall, that under the
isomotphism $\rep{W}+\rep{E}\sim V^{-2}$ vectors $u_s, k_s,
f_1,f_p,f_{p+1}, f_{2p}$ correspond to the linearly independent
functions \eqref{Yeigenfunctions-u(x)}--\eqref{Yeigenfunctions-f(x)}
in~$V^{-2}$, and therefore the vectors are also linearly independent.
In particular, vectors $u_s, f_1,f_p,f_{p+1}, f_{2p}$ form a basis in
$\rep{W}$ and therefore images of the vectors $k_s$ under the
canonical projection to $\rep{E}=(\rep{W}+\rep{E})/\rep{W}$ form a
basis in $\rep{E}$. We let abusing notations $k_s$ denote these
images. We recall that the isomorphism $\rep{E}\sim \rep{M}$ maps
vectors $k_s$ (images under the canonical projections of
$k_s\in\rep{Z}$) to $f_s$ (images under the canonical projections of
$f_s\in\rep{Z}$), and therefore $f_s$ with $s=2\dots p-1, p+2 \dots
2p-1$ are linearly independent. Thus, the linear independence of all
vectors $u,f,k$ is established.
\end{proof}

\section{$PSL(2,\oZ)$ action in $\rep{Z}$\label{sec:psl2z}}
In this section we define operators $\modS$ and $\ribbon$ and prove
that they satisfy $PSL(2,\oZ)$ relations.  Conjugations with operators
$\modS$ and $\ribbon$ give automorphisms $\sigma$ and $\tau_+$
respectively.
\subsection{$\sigma$\label{sec:Strans}}
 We define the $\modS$-operator that maps the $X$-basis to the $Y$-basis by
formulas~\eqref{S-action}.
\begin{prop}
$\modS$ satisfies relations
\begin{gather}
\label{SXSinZ}\modS X\modS^{-1}=Y^{-1},\\
\label{SYSinZ} \modS Y\modS^{-1}=XT^2,\\
 \label{STSinZ} \modS T\modS ^{-1}=T,\\
\label{S^2inZ}   \modS ^2=\q T^{-1}.
\end{gather}
\end{prop}
\begin{proof}
\begin{itemize}
\item \eqref{SXSinZ} follows from the definition of $\modS$.
\item \eqref{S^2inZ} follows from a direct calculation of
  $T\modS^2$-action in the $X$-basis.  We give a detailed calculation of
  $T\modS^2e_s$. The calculation of $T\modS^2w_s$ and $T\modS^2m_s$ is
  similar and is omitted. We check that $T\modS^2e_s=\q e_s.$ We begin
  with
\begin{gather}
\notag \modS^2e_s=\modS f_s=\modS
\left(\sum_{r=1}^{2p}\fdec_{r,s}^{(w)}w_r+\fdec_{p,s}^{(e)}e_p
+\fdec_{2p,s}^{(e)}e_{2p}+\left(\sum_{r=2}^{p-1}
+\sum_{r=p+2}^{2p-1}\right)\fdec_{r,s}^{(m)}m_r\right)=\\
\label{S^2e_s}
\begin{split}=\fdec_{1,s}^{(w)}u_1+\fdec_{p,s}^{(w)}u_p
+\fdec_{p+1,s}^{(w)}u_{p+1}&+\fdec_{2p,s}^{(w)}u_{2p}
+\fdec_{p,s}^{(e)}f_p+\fdec_{2p,s}^{(e)}f_{2p}+\\
 +\left(\sum_{r=2}^{p-1}
+\sum_{r=p+2}^{2p-1}\right)&(\fdec_{r,s}^{(w)}u_r+\fdec_{r,s}^{(m)}k_r).\end{split}
\end{gather}
Then we calculate coefficients in front of $e_j$ and $w_j$ in
\eqref{S^2e_s} using \eqref{Yeigenvectors-u}-\eqref{f-coefficients}.
This calculation is cumbersome and is given in
Appendix~\ref{appendixA}. The result of the calculation is
\begin{gather}
\begin{split}
\modS^2e_s=-\frac{\q^{s+1}}{[s]}e_s+\q\frac{[s-1]}{[s]}e_{2p-s}
  +&\frac{2(\q^2-1)}{(\q^s-\q^{-s})^2}(w_s-w_{2p-s}),\quad s\ne p,2p,\\
\modS^2 e_p=e_p,&\quad \modS^2 e_{2p}=e_{2p}.
\end{split}
\end{gather}
A simple calculation using \eqref{Tact-w_p}-\eqref{Tact-e_s} gives
$$T\modS^2e_s=\q e_s,\quad s=1\dots 2p.$$
\item \eqref{STSinZ} is checked as follows
$\modS ^2=\q T^{-1}\Rightarrow \modS
T=T\modS(=\q\modS^{-1})\Rightarrow \modS T\modS ^{-1}=T.$
 \item \eqref{SYSinZ} is checked as follows
$\modS X\modS^{-1}\stackrel{\eqref{SXSinZ}}{=}Y^{-1}\Rightarrow \modS
X^{-1}\modS^{-1}=Y\Rightarrow \modS Y\modS^{-1}=\modS^2
X^{-1}\modS^{-2}\stackrel{\eqref{S^2inZ}}{=}T^{-1}X^{-1}T
\stackrel{\eqref{dha1}}{=}XT^2$.
\end{itemize}
 \end{proof}

\subsection{$\tau_+$\label{sec:tau}}
The automorphism $\tau_+$ can be realized as a conjugation with the element
$\ribbon\in\rep{Z}$ given by~\eqref{ribbon}.
\begin{prop}
For $\ribbon$ given by~\eqref{ribbon}, the operator
\begin{equation*}
 \tau_+(x)=\ribbon^{-1} x \ribbon,\quad\forall x\in\daha
\end{equation*}
satisfies relations \eqref{taupl}.
\end{prop}
\begin{proof}
  A direct calculation.
\end{proof}

\begin{prop}
The map
  \begin{equation*}
 {\scriptstyle\left(
\begin{array}{cc}
0&1\\
-1&0
\end{array} \right)}\rightarrow \modS,\qquad
 {\scriptstyle\left(
\begin{array}{cc}
1&1\\
0&1
\end{array} \right)}\rightarrow \ribbon\cdot
\end{equation*}
gives a $PSL(2,\oZ)$ action in $\rep{Z}$.
\end{prop}
\begin{proof}
  Relations \eqref{SXSinZ}--\eqref{S^2inZ} and \eqref{taupl} are
  sufficient~\cite{[Cherednik-book]} to check the $PSL(2,\oZ)$
  relations.
\end{proof}

\section{eigenspace of $T$ with eigenvalue $\q$\label{sec:symmetrization}}
In this section, we describe the representation of the symmetrized
DAHA. In section~\ref{Proof:main}, we prove that it is isomorphic to
the centre of $\UresSL2$.

We let $\rep{T_\q}$ denote the eigenspace of~$T$ with the eigenvalue
$\q$. It is $3p-1$  dimensional. Operators $X$ and $Y$ have no
well-defined restriction to $\rep{T_\q}$ but "symmetrized"
operators $C=-(X+X^{-1})$ and $H=-(Y+Y^{-1})$ have. Indeed, for a
given $\mathbf{a}\in\rep{T_\q}$, we have
\begin{multline}
  T(X+X^{-1})\mathbf{a}\stackrel{\eqref{dha1}}{=}(X^{-1}T^{-1}
+TX^{-1})\mathbf{a}=(\q^{-1}+T)X^{-1}\mathbf{a}
\stackrel{\eqref{dha2}}{=}\\
=(\q+T^{-1})X^{-1}\mathbf{a}=(\q X^{-1}+XT)\mathbf{a}
=\q(X+X^{-1})\mathbf{a}.
\end{multline}
Thus, $(X+X^{-1})\mathbf{a}\in \rep{T_\q}$ and a similar calculation
shows that $H$ has the well-defined restriction to~$\rep{T_\q}$ as
well.

\subsection{$C$-basis}
The eigenvectors of $C=-(X+X^{-1})$ are
\begin{gather*}
\idem_{0}=e_{p},\quad\idem_p=e_{2p},\quad \idem_s=e_{p+s}+e_{p-s},\\
\nilp_1^+=\frac{1}{\q-\q^{-1}}m_{p-1},\quad\nilp_1^-=\frac{1}{\q-\q^{-1}}
\left(w_{p+1}-w_{p-1}-m_{p-1}\right), \\
\nilp_s^+=\frac{[s]}{\q-\q^{-1}}\left(m_{p-s}+m_{p+s}\right),\quad \nilp_s^-
=\frac{[s]}{\q-\q^{-1}}\left(w_{p+s}-w_{p-s}-m_{p-s}-m_{p+s}\right),\\
\nilp_{p-1}^+=\frac{1}{\q-\q^{-1}}\left(m_{2p-1}-w_1\right),\quad{}\nilp_{p-1}^-
=\frac{1}{\q-\q^{-1}}(w_{2p-1}-m_{2p-1}),\\
\nilp_s=\nilp_s^++\nilp_s^-=\frac{[s]}{\q-\q^{-1}}(w_{p+s}-w_{p-s}),\quad
s=1\dots p-1.
\end{gather*}
The action of $C$ on them folows from \eqref{X-jordanaction}
\begin{gather}
C\idem_0=\mu_0\idem_0,\quad C\idem_p=\mu_p\idem_p,\\
C\idem_s=\mu_s\idem_s+(\q-\q^{-1})^2\nilp_s,\quad  \qquad s=1\dots p-1\\
C\nilp_s^\pm=\mu_s\nilp^\pm_s,\qquad s=1\dots p-1,
\end{gather}
where
\begin{equation}\label{mu}
 \mu_s= \q^s+\q^{-s},\qquad0\leq s\leq p.
\end{equation}
The multiplication in $\rep{T_\q}$ is induced by \eqref{multipl-I}
\begin{equation}
\label{multipl-II}
\idem_r\nilp_s^{\pm}=\delta_{r,s}\nilp_s^{\pm},\qquad\idem_r\idem_s=\delta_{r,s}\idem_s,\qquad
\nilp_r^{\pm}\nilp_s^{\pm}=0.
\end{equation}

\subsection{$H$-basis}
The eigenvectors of $H=-(Y+Y^{-1})$ are
\begin{gather*}
\fidem_{0}=f_{p},\quad \fidem_p=f_0,\quad \fidem_s=f_{p+s}+f_{p-s},\\
\unilp_1^+=\frac{1}{\q-\q^{-1}}k_{p-1},\quad\unilp_1^-=\frac{1}{\q-\q^{-1}}
\left(u_{p+1}-u_{p-1}-k_{p-1}\right), \\
\unilp_s^+=\frac{[s]}{\q-\q^{-1}}\left(k_{p-s}+k_{p+s}\right),\quad \unilp_s^-
=\frac{[s]}{\q-\q^{-1}}(u_{p+s}-u_{p-s}-k_{p-s}-k_{p+s}),\\
\unilp_{p-1}^+=\frac{1}{\q-\q^{-1}}\left(k_{2p-1}-u_1\right),\quad{}
\unilp_{p-1}^-=\frac{1}{\q-\q^{-1}}(u_{2p-1}-k_{2p-1}),\\
\unilp_s=\unilp_s^++\unilp_s^-=\frac{[s]}{\q-\q^{-1}}(u_{p+s}-u_{p-s}),\quad
s=1\dots p-1.
\end{gather*}
The action of $H$ on them folows from \eqref{Y-jordanbasis}
\begin{gather}
H\fidem_0=\mu_0\fidem_0,\quad H\fidem_p=\mu_p\fidem_p, \\
H\fidem_s=\mu_s\fidem_s+(\q-\q^{-1})^2\unilp_s,\quad s=1\dots p-1,\\
H\unilp_s^\pm=\mu_s\unilp^\pm_s, \qquad s=1\dots p-1,
\end{gather}
where eigenvalues are given by~\eqref{mu}.

\subsection{$SL(2,\oZ)$ action}
Operators $\modS$ and $\ribbon$ have well-defined restrictions to
$\rep{T_\q}$. This endows $\rep{T_\q}$ with a representation of
$SL(2,\oZ)$. In more detail, $\modS$-operator in $\rep{T_\q}$
satisfies
\begin{equation}\label{SinTq}
  \begin{split}
   \modS\idem_s&=\fidem_s,\quad s=0\dots p,\\
   \modS\nilp^\pm_s&=\unilp^\pm_s,\quad s=1\dots{}p-1
  \end{split}
\end{equation}
and because $T=\q$ in $\rep{T_\q}$, we have $\modS^2=1$.  We note also
that in $\rep{T_\q}$ relations \ref{SXSinZ} and \ref{SYSinZ} lead to
\begin{equation}
  \modS C \modS^{-1}=H.
\end{equation}

\section{Proof of Theorem~\ref{Thm:main}}
\label{Proof:main} We note that $\rep{T_\q}$ and the center of
$\UresSL2$ from \cite{[FGST-mod]} have the same dimension equal
to~$3p-1$. Then we identify $C$-basis ($H$-basis) in $\rep{T_\q}$ with
the Radford images (Drinfeld images) of $q$-characters of irreducible
representations
 \begin{equation} \label{notations-2}
   \begin{split}
\hat{\radmap}^+(s)&=\omega_s\nilp_s^+,\quad
\hat\radmap^-(s)=\omega_{p-s}\nilp_{p-s}^-,\quad s=1\dots p-1,\\
\hat\radmap^+(p)&=p\sqrt{2p}\idem_p,\quad \hat\radmap^-(p)=(-1)^{p+1}p\sqrt{2p}\idem_0,\\
\cchi^+(s)&=\omega_s\unilp_s^+,\quad
\cchi^-(s)=\omega_{p-s}\unilp_{p-s}^-,\quad s=1\dots p-1,\\
\cchi^+(p)&=p\sqrt{2p}\fidem_p,\quad \cchi^-(p)=(-1)^{p+1}p\sqrt{2p}\fidem_0.
   \end{split}
 \end{equation}
 This identification establishes an isomorphism between $\rep{T_\q}$
 and the center of $\UresSL2$ as associative commutative algebras.

 Under the identification \eqref{notations-2}, $\rep{T_\q}$ coincides
 with the center of $\UresSL2$ as the representation of $SL(2,\oZ)$.
 In particular, the relations
 $\modS(\cchi^\pm(s))=\hat\radmap^\pm(s)$ for $s=0\dots p$ in the
 center are parallel to the relations~\eqref{SinTq} in~$\rep{T_\q}$.
 The Gaussian element $\ribbon$ in notations of~\cite{[FGST-mod]}
\begin{equation}
  \label{v}
\ribbon=\sum_{s=0}^{p}(-1)^{s+1}\q^{-\frac{1}{2}(s^2-1)}\idem_s
+\sum_{s=1}^{p-1}(-1)^{p}\q^{-\frac{1}{2}(s^2-1)}\frac{\q^s-\q^{-s}}{\sqrt{2p}}\vvarphi(s),
\end{equation}
where
$\vvarphi(s)=\frac{p-s}{p}\hat\radmap^+(s)-\frac{s}{p}\hat\radmap^-(p-s)$
for $1\le s\le p-1$ coincides with the ribbon element of $\UresSL2$.



\section{Discussion}
We identified the representation of DAHA that gives the Verlinde
algebra of $(1,p)$ logarithmic conformal field models. The center of
$\UresSL2$ coincides with the symmetrization of $\rep{Z}$ and
$C=-(X+X^{-1})$ coincides with the $\UresSL2$ Casimir element.
Probably the whole representation $\rep{Z}$ can be realized in
$\UresSL2$ such that $X$ would be realized by a multiplication with a
$\UresSL2$ element.

Another interesting direction of investigations is to find a
realization of $\daha$ on $(1,p)$ logarithmic conformal field
model conformal blocks. This can also be useful in boundary
conformal field theories.  The Ishibashi and Cardy boundary states
can probably be identified with eigenvectors of operators
$C=-(X+X^{-1})$ and $H=-(Y+Y^{-1})$ respectively.
\subsubsection*{Acknowledgments}
We are grateful to T.~Suzuki, M.~Kasatani and T.~Kuwabara for many
valuble discussions and A.M.~Semikhatov for the useful discussions and
comments on a presentation of our results.  The work of GM was
supported by the RFBR Grant 07-01-00523. The work of IYuT was
supported in part by LSS-4401.2006.2 grant, the RFBR Grant 05-02-17217
and the ``Dynasty'' foundation.

\appendix
\section{proof of $\modS^2=\q T^{-1}$}
\label{appendixA} We calculate coefficient in front of $e_j$ in
\ref{thefactorofe_j} and coefficient in front of $w_j$ in
\ref{thefactorofw_j}.
\subsection{The coefficient in front of $e_j$\label{thefactorofe_j}}
The substitution of \eqref{Yeigenvectors-u},
\eqref{Yeigenvectors-k}, \eqref{Yeigenvectors-f} in \eqref{S^2e_s}
gives the coefficient in front of $e_j$
\begin{multline}
\label{e_j-factor}
\fdec_{1,s}^{(w)}\udec_{j,1}^{(e)}+\fdec_{p,s}^{(w)}\udec_{j,p}^{(e)}
+\fdec_{p+1,s}^{(w)}\udec_{j,p+1}^{(e)}+\fdec_{2p,s}^{(w)}\udec_{j,2p}^{(e)}
+\fdec_{p,s}^{(e)}\fdec_{j,p}^{(e)}+\fdec_{2p,s}^{(e)}\fdec_{j,2p}^{(e)}+\\
+ \left(\sum_{r=2}^{p-1}+\sum_{r=p+2}^{2p-1}\right)
(\underbrace{\fdec_{r,s}^{(w)}\udec_{j,r}^{(e)}+\fdec_{r,s}^{(m)}\kdec_{j,r}^{(e)}}_{A}),
 \end{multline}
where all numbers $u$,  $k$, $f$  are given in
\eqref{u-coefficients}, \eqref{k-coefficients},
\eqref{f-coefficients}. A simplification of the underbraced
expression gives
\begin{description}
\item[for $j\ne1,p,p+1,2p$]
  \begin{multline*}
A=\frac{(-1)^{s+j}\q^2}{2p}\Big(\q^{r-1}[s,r-1][r,j]-[s,r-1][r,j-1]\Big)+\\
  +\frac{(-1)^{s+j}\q^s}{2p}\Big(\q[s,r][r,j-1]-\q^r[s,r][r,j]\Big),
  \end{multline*}
\item[for $j=1$]
  \begin{equation*}
A=\frac{(-1)^{s+1}\q^2}{2p}\Big(\q^{r-1}[s,r-1][r,1]-[s,r-1]\Big)
   +\frac{(-1)^{s+1}\q^s}{2p}\Big(-\q\{s,r\}-\q^r[s,r][r,1]\Big),
  \end{equation*}
\item[for $j=p+1$]
  \begin{multline*}
A=\frac{(-1)^{s+p+1}\q^2}{2p}\Big(\q^{r-1}[s,r-1][r,p+1]-[s+p,r-1]\Big)+\\
  +\frac{(-1)^{s+p+1}\q^s}{2p}\Big(\q\{s+p,r\}-\q^r[s,r][r,p+1]\Big),
\end{multline*}
\item[for $j=p$]
  \begin{multline*}
A=\frac{(-1)^{s+p}\q^2}{2p}\Big(\q^{r-1}[s+p,r-1]-[s,r-1][r,p-1]\Big)+\\
  +\frac{(-1)^{s+p}\q^s}{2p}\Big(\q[s,r][r,p-1]-\q^r\{s+p,r\}\Big),
\end{multline*}
\item[for $j=2p$]
  \begin{equation*}
A=\frac{(-1)^{s}\q^2}{2p}\Big(\q^{r-1}[s,r-1]-[s,r-1][r,2p-1]\Big)
  +\frac{(-1)^{s}\q^s}{2p}\Big(\q[s,r][r,2p-1]+\q^r\{s,r\}\Big).
  \end{equation*}
\end{description}
Then the simplification of \eqref{e_j-factor} gives coefficients
 in \eqref{S^2e_s} in front of $e_s$. Explicitly, the summation in $r$ of different terms in $A$
is given by
 \begin{gather}
\left(\sum\limits_{r=2}^{p-1}
  +\sum\limits_{r=p+2}^{2p-1}\right)q^{r-1}[s,r-1][r,j]
   = \Bigl(p-\frac{1}{4}\bigl((s+j+1\,\mbox{\rm mod}\,2p)+(s-j+1\,\mbox{\rm mod}\,2p)+\notag\\
   +(s+j-1\,\mbox{\rm mod}\,
  2p)+(s-j-1\,\mbox{\rm mod}\, 2p)\bigr)\Bigr)\Big(1+(-1)^{s+j}\Big)
+p\frac{\{j\}}{[j]}(\delta_{s+j,2p}-\delta_{s-j,0}),\notag\\
\kern-30pt\left(\sum\limits_{r=2}^{p-1}+\sum\limits_{r=p+2}^{2p-1}\right)[s,r-1][r,j]
=\Big(1+(-1)^{s+j+1}\Big)\Big(p-\frac{1}{2}\Big((s+j\,\mbox{\rm mod}\,2p)+\notag\\
        \kern300pt+(s-j\,\mbox{\rm mod}\,2p)\Big)\Big)\notag\\
\left(\sum\limits_{r=2}^{p-1}+\sum\limits_{r=p+2}^{2p-1}\right)[s,r][r,j]
=[s]((-1)^{s+j}-1),\quad j\ne p,2p\notag\\
\left(\sum\limits_{r=2}^{p-1}+\sum\limits_{r=p+2}^{2p-1}\right)\q^r[s,r][r,j]
=\frac{2p(\delta_{s,2p-j}-\delta_{s,j})}{\q^j-\q^{-j}}-\q[s](1+(-1)^{j+s}),\quad j\ne p,2p,\notag\\
\left(\sum\limits_{r=1}^{p-1}+\sum\limits_{r=p+1}^{2p-1}\right)[s,r]
  =(1-(-1)^s)\bigl(p-(s\,\mbox{\rm mod}\,2p)\bigr),\notag\\
\sum\limits_{r=1}^{2p}\{s,r\}=0,\notag\\
\left(\sum\limits_{r=1}^{p-1}+\sum\limits_{r=p+1}^{2p-1}\right)\q^r[s,r]
=(1+(-1)^s)\left(p-\frac{((s+1)\,\mbox{\rm mod}\,2p)+((s-1)\,\mbox{\rm mod}\,2p)}{2}\right),\notag\\
\sum\limits_{r=1}^{2p}\q^r\{s,r\}=(1+(-1)^s)\left(\frac{((s-1)\,\mbox{\rm mod}\,2p)
             -((s+1)\,\mbox{\rm mod}\,2p)}{2}\right).\notag
\end{gather}

\subsection{The coefficient in front of $w_j$\label{thefactorofw_j}}
The substitution of \eqref{Yeigenvectors-u},
\eqref{Yeigenvectors-k}, \eqref{Yeigenvectors-f} in \eqref{S^2e_s}
gives the coefficient in front of $w_j$
\begin{multline}\label{w_j-factor}
\fdec_{1,s}^{(w)}\udec_{j,1}^{(w)}+\fdec_{p,s}^{(w)}\udec_{j,p}^{(w)}
+\fdec_{p+1,s}^{(w)}\udec_{j,p+1}^{(w)}+\fdec_{2p,s}^{(w)}\udec_{j,2p}^{(w)}
+\fdec_{p,s}^{(e)}\fdec_{j,p}^{(w)}+\fdec_{2p,s}^{(e)}\fdec_{j,2p}^{(w)}+\\
+ \left(\sum_{r=2}^{p-1}+\sum_{r=p+2}^{2p-1}\right)
(\underbrace{\fdec_{r,s}^{(w)}\udec_{j,r}^{(w)}+\fdec_{r,s}^{(m)}\kdec_{j,r}^{(w)}}_{A}),
\end{multline}
where all numbers $u$,  $k$, $f$  are given in
\eqref{u-coefficients}, \eqref{k-coefficients},
\eqref{f-coefficients}. A simplification of the underbraced
expression gives
\begin{multline}\label{theA}
  A=\frac{\q^2(-1)^{s+j}\{1,j-1\}}{p(\q^{j-1}+\q^{-j+1})}[s,r-1][r-1,j-1]
       -\frac{(-1)^{s+j}\{1,j\}\q^{2}}{p(\q^j+\q^{-j})}\q^{r-1}[s,r-1][r-1,j]+\\
+\frac{\q^s(-1)^{s+j}}{2p}\Big(\q^r\{s,r\}\{r,j\}-\q\{s,r\}\{r,j-1\}
  +\q^r[s,r][r,j]\{1,j\}-\q [s,r][r,j-1]\{1,j-1\}\Big)
\end{multline}
Then the simplification of \eqref{w_j-factor} gives coefficients
 in \eqref{S^2e_s} in front of $w_s$. Explicitly, the summation in $r$ of different terms
in~\eqref{theA} is given by
 \begin{gather*}
 \left(\sum\limits_{r=2}^{p-1}+\sum\limits_{r=p+2}^{2p-1}\right)[s,r-1][r-1,j]
           =[s](1-(-1)^{s+j}),\quad{}j\ne{}p,2p,\notag\\
 \left(\sum\limits_{r=2}^{p-1}+\sum\limits_{r=p+2}^{2p-1}\right)\q^{r-1}[s,r-1][r-1,j]
 =\frac{2p(\delta_{s,2p-j}-\delta_{s,j})}{\q^j-\q^{-j}}+\notag\\
                    \kern250pt+ \q^{-1}[s](1+(-1)^{j+s}),\quad j\ne p,2p,\notag\\
 \left(\sum\limits_{r=2}^{p-1}+\sum\limits_{r=p+2}^{2p-1}\right)\q^r\{s,r\}\{r,j\}
 =\Big(\frac{p(\delta_{s+j,2p}+\delta_{s-j,0})-2}{\q^j+\q^{-j}}-\q\{s\}\Big)\{1,j\}(1+(-1)^{j+s}),\notag\\
 \left(\sum\limits_{r=2}^{p-1}+\sum\limits_{r=p+2}^{2p-1}\right)\{s,r\}\{r,j\}
      =\{s,1\}\{1,j\}((-1)^{s+j}-1),\notag\\
 \left(\sum\limits_{r=2}^{p-1}+\sum\limits_{r=p+2}^{2p-1}\right)\q^r[s,r][r,j]
 =\frac{2p(\delta_{s,2p-j}-\delta_{s,j})}{\q^j-\q^{-j}}-\q[s](1+(-1)^{j+s}),\quad j\ne p,2p,\notag\\
 \left(\sum\limits_{r=2}^{p-1}+\sum\limits_{r=p+2}^{2p-1}\right)[s,r][r,j]
        =[s]((-1)^{s+j}-1),\quad j\ne p,2p.\notag
 \end{gather*}

\end{document}